\documentclass[10pt,a4paper]{amsart}
\usepackage{amsmath}
  \usepackage{paralist}
  \usepackage{graphics} %% add this and next lines if pictures should be in esp format
  \usepackage{epsfig} %For pictures: screened artwork should be set up with an 85 or 100 line screen
 \usepackage[colorlinks=true]{hyperref}
\hypersetup{urlcolor=blue, citecolor=red}

\newtheorem{theorem}{Theorem}[section]

\newtheorem{lemma}[theorem]{Lemma}

\theoremstyle{definition}

\newtheorem{remark}{Remark}

\numberwithin{equation}{section}

\newcommand{\ep}{\varepsilon}
\newcommand{\R}{\mathbb{R}}
\newcommand{\am}{\alpha _-}
\newcommand{\ap}{\alpha _+}
\newcommand{\apm}{\alpha _ \pm}

\title[Validity of expansions in Allen-Cahn and FitzHugh-Nagumo]
      {On the validity of formal asymptotic
expansions in Allen-Cahn equation and FitzHugh-Nagumo system with
generic initial data}

\author[M. Alfaro and H. Matano]{}

\subjclass{35B25, 35C20, 35R35} \keywords{Singular perturbation,
asymptotic expansion, front profile, Allen-Cahn equation,
FitzHugh-Nagumo system, reaction-diffusion system}

 \email{malfaro@math.univ-montp2.fr}
 \email{matano@ms.u-tokyo.ac.jp}

\thanks{The first author is supported by the French Agence
Nationale de la Recherche within the project IDEE
(ANR-2010-0112-01), and the second author by KAKENHI (23244017)}

\begin{document}
\maketitle

\centerline{\scshape Matthieu Alfaro }
\medskip
{\footnotesize
 \centerline{Univ. Montpellier 2, I3M, UMR CNRS 5149,}
 \centerline{CC051,
 Place Eug\`ene Bataillon, 34095 Montpellier Cedex 5, France,}
}

\medskip

\centerline{\scshape Hiroshi Matano}
\medskip
{\footnotesize \centerline{Univ. of Tokyo, Graduate School of
Mathematical Sciences,} \centerline{ 3-8-1 Komaba, Tokyo 153-8914,
Japan.} }

\begin{abstract}
Formal asymptotic expansions have long been used to study the
singularly perturbed Allen-Cahn type equations and
reaction-diffusion systems, including in particular the
FitzHugh-Nagumo system. Despite their successful role, it has been
largely unclear whether or not such expansions really represent
the actual profile of solutions with rather general initial data.
By combining our earlier result and known properties of eternal
solutions of the Allen-Cahn equation, we prove validity of the
principal term of the formal expansions for a large class of
solutions.
\end{abstract}

%%%%%%%%%%%%%%%%%%%%%%%%%%%%%%%
%%%%%%%%%%%%%%%%%%%%%%%%%%%%%%%
\section{Introduction}\label{s:intro}

In this paper, we study the behavior of solution $u^\ep$ of an
Allen-Cahn type equation of the form
\[
 ({\rm P}^\ep) \quad\begin{cases}
 u_t=\Delta u+\displaystyle \frac 1{\ep ^2} (f(u)-\ep g^\ep(x,t,u))
 &\text{in }\Omega \times (0,\infty)\vspace{2pt}\\
\displaystyle \frac{\partial u}{\partial \nu} = 0
&\text{on }\partial \Omega \times (0,\infty)\vspace{2pt}\\
 u(x,0)=u_0(x) &\text{in }\Omega\,,
 \end{cases}
\]
and also that of a reaction-diffusion system of the form
\[
 ({\rm RD} ^\ep) \quad
 \begin{cases}
 u_t=\Delta u+\displaystyle \frac 1{\ep ^2}\,(f(u)+\ep f_1(u,v)+
 \ep^2 f_2 ^\ep (u,v)) \quad&\textrm{in}\ \
 \Omega\times(0,\infty)\vspace{2pt}\\
 v_t=D \Delta v +h(u,v) \quad&\textrm{in} \ \
 \Omega\times(0,\infty)\vspace{2pt}\\
\displaystyle \frac{\partial u}{\partial \nu}
=\displaystyle{\frac{\partial v}{\partial\nu}}
 = 0 \quad &\textrm{on}\ \
 \partial\Omega\times(0,\infty) \vspace{2pt}\\
 u(x,0)=u_0(x)\quad &\textrm{in}\ \ \Omega\vspace{2pt}\\
 v(x,0)=v_0(x)\quad &\textrm{in}\ \ \Omega\,,
 \end{cases}\hspace{30pt}
\]
where  $f(u), g^\ep(x,t,u), f_1(u,v)$ and $f_2^\ep(u,v)$ satisfy
the conditions to be specified later, and $\ep$ is a positive
parameter. A typical example of $({\rm RD}^\ep)$ is the
FitzHugh-Nagumo system:
\[
 ({\rm FHN}^\ep) \quad
 \begin{cases}
 \,u_t=\Delta u+\displaystyle \frac 1{\ep^2}\,
 (f(u)-\ep f_1 (u)-\ep v) \quad&\textrm{in}\ \
 \Omega\times(0,\infty)\vspace{2pt}\\
 \,v_t=D \Delta v +\alpha u -\beta v \quad&\textrm{in} \ \
 \Omega\times(0,\infty)\vspace{2pt}\\
 \,\displaystyle\frac{\partial u}{\partial \nu}
 =\displaystyle{\frac{\partial v}{\partial\nu}}
 = 0 \quad &\textrm{on}\ \
 \partial\Omega\times(0,\infty) \vspace{2pt}\\
 \,u(x,0)=u_0(x)\quad &\textrm{in}\ \ \Omega\\
 \,v(x,0)=v_0(x)\quad &\textrm{in}\ \ \Omega\,.
 \end{cases}\hspace{30pt}
\]

It is well-known that the solution $u^\ep(x,t)$ of the above
systems develops a steep transition layer, which converges to a
``sharp interface" as $\ep\to 0$. To study such a sharp interface
limit, formal asymptotic expansions of $u^\ep$ are commonly used
to discover, formally, the law of motion of the limit interface.
Then, based on these expansions, one can construct sub- and
super-solutions or use some approximation argument to prove the
convergence of the transition layer --- or the front --- to the
sharp interface, thereby establishing rigorously that the limit
motion law agrees with what is anticipated from the formal
asymptotics.

However, this standard approach only tells us that the transition
layer of $u^\ep$ is confined within a relatively narrow zone ---
of thickness $o(1)$ or sometimes even ${\mathcal O}(\ep)$---
around the limit interface, but it does not say much about whether
or not the actual transition layer really possesses a robust
profile that matches the formal asymptotics. Known answers to this
question are mainly concerned with solutions whose initial data
already has a well-developed transition layer. The case of more
general solutions has largely been unexplored. Our goal is to
provide an affirmative answer in this direction: we shall prove
that, for $\ep$ sufficiently small, the solution $u^\ep$ --- with
rather general initial data --- of both $($P$^\ep)$ and $({\rm
RD}^\ep)$ possesses a profile that agrees with the principal term
of the formal expansion.

%%%%%%%%%%%%%%%%%%%%%%%%%%%
\subsection{Notation and assumptions}
\label{ss:assumptions}

The notation and assumptions stated below strictly follow those in
\cite{A-Hil-Mat}.

In problems $({\rm P}^\ep)$ and $({\rm RD}^\ep)$ above, $\Omega$
is a smooth bounded domain in $\R^N$ ($N\geq 2$) and $\nu$ denotes
the outward unit normal vector to $\partial \Omega$. The
nonlinearity is given by $f(u):=-W'(u)$, where $W(u)$ is a
double-well potential with equal well-depth, taking its global
minimum value at $u=\apm$. More precisely we assume that $f$ is
$C^2$ and has exactly three zeros $\am<a<\ap$ such that
\begin{equation}\label{der-f}
f'(\apm)<0\,, \quad f'(a)>0\,\quad\ \hbox{(bistable
nonlinearity)},
\end{equation}
and that
\begin{equation}\label{int-f}
\int _ {\am} ^ {\ap} f(u)\,du=0\,\quad \hbox{ (balanced case)}.
\end{equation}

In the Allen-Cahn equation $({\rm P}^\ep)$ we allow the balance of
the two stable zeros $\am$ and $\ap$ to be slightly broken by the
function  $-\ep g^\ep(x,t,u)$ defined on
$\overline{\Omega}\times[0,\infty)\times\R$. We assume that
$g^\ep$ is $C^2$ in $x$ and $C^1$ in $t,\,u$, and that, for any
$T>0$ there exist $\vartheta \in (0,1)$ and $C>0$ such that, for
all $(x,t,u) \in\overline{\Omega}\times[0,T]\times \R$,
\begin{equation*}
 \ep |\Delta_{x} g^\ep (x,t,u)| + \ep| g^\ep _t(x,t,u)|+
 |g^\ep _u (x,t,u)| \leq C\,,
\end{equation*}
\begin{equation*}
 \Vert g^\ep(\cdot,\cdot,u)
 \Vert_{C^{1+\vartheta,\frac{1+\vartheta}{2}}(\bar \Omega \times[0,T])}
 \leq C\,.
\end{equation*}
Moreover, we assume that there exists a function $g(x,t,u)$ and a
constant, which we denote again by $C$, such that
\begin{equation}\label{def:g}
 |g^\ep (x,t,u)-g(x,t,u)| \leq C \ep\,,
\end{equation}
for all small $\ep>0$. In \cite{A-Hil-Mat}, we also assumed
$\frac{\partial  g^\ep}{\partial \nu} =0$ on $\partial \Omega
\times [0,\infty)\times \mathbb R$, though this last condition is
only for technical simplicity.  Note that these conditions, except
the last one, are automatically satisfied if $g^\ep$ is smooth and
independent of $\ep$.

In the reaction-diffusion system $({\rm RD}^\ep)$ we assume that
$f_1(u,v)$, $f_2^\ep(u,v)$ are $C^2$ functions and that $f_2^\ep$,
along with its derivatives, remain bounded as $\ep \to 0$. We also
assume that $D>0$ and that $h(u,v)$ is a $C^2$ function such that,
for any constants $L,M>0$, there exists a constant $M_1\geq M$
such that $ h(u,-M_1)\geq 0\geq h(u,M_1)$ for $ |u|\leq L$. These
conditions enable us to construct a family of invariant rectangles
mentioned in Remark \ref{rem:rectangles}.

In the FitzHugh-Nagumo system, we assume that $f_1 \in C^2(\R)$
and that $\alpha$ and $\beta$ are given positive constants so that
$({\rm FHN}^\ep)$ becomes a special case of the reaction-diffusion
system $({\rm RD} ^\ep)$.

To complete the picture we need to specify conditions on the
initial data. We assume that $u_0$ and $v_0$ belong to
$C^2(\bar\Omega)$. We define the ``initial interface" $\Gamma _0$
by
\begin{equation}
\Gamma_0:=\{x\in\Omega:\,u_0(x)=a\}\,,
\end{equation}
and assume that $\Gamma_0$ is a $C^{3+\vartheta}$ hypersurface
($0<\vartheta<1$) without boundary such that
\begin{equation*}
\Gamma_0 \subset\subset \Omega \quad \mbox { and } \quad \nabla
u_0(x)\cdot n(x,0) \neq 0\quad\text{if $x\in\Gamma_0\,,$}
\end{equation*}
\begin{equation*}
u_0>a \quad \text { in }\Omega ^+ _0\,,\quad u_0<a \quad \text {
in }\Omega ^- _0\,,
\end{equation*}
where $\Omega ^- _0$ denotes the region enclosed by $\Gamma _0$
and $\Omega ^+ _0$ the region enclosed between $\partial \Omega$
and $\Gamma_0$, and $n(x,0)$ denotes the outward unit normal
vector at $x\in \Gamma_0=\partial \Omega^-_0$. Let us emphasize
that we do not assume that the initial data $u_0$ of $u^\ep$
already has well-developed transition layers depending on $\ep$,
in which case the validity of the formal expansions is more or
less known (see subsection \ref{ss:singular} for more details).

\begin{remark}[Time-global smooth solutions]\label{rem:rectangles}
Under the above assumptions, it is classical that $({\rm P}^\ep)$
has a uniformly bounded smooth solution $u^\ep$ that exists for
all $t\geq 0$. As for $({\rm RD} ^\ep)$, the same can be shown for
$\ep >0$ small enough, by using the method of invariant rectangles
(see e.g. \cite{A-Hil-Mat} for details).
\end{remark}

%%%%%%%%%%%%%%%%%%%%%%%%%%%%%%%%%%
\subsection{Known results for the singular limit}
\label{ss:singular}

We present here a brief overview of known results. Heuristically,
in the very early stage, the diffusion term  is negligible
compared with the reaction term. Hence, in view of the profile of
$f$, the value of $u^{\ep}$ quickly becomes close to either $\ap$
or $\am$ in most part of $\Omega$, creating a steep interface
(transition layers) between the regions $\{u^\ep\approx \am\}$ and
$\{u^\ep\approx \ap\}$ ({\it Generation of interface}). Once the
balance between diffusion and reaction near the transition layers
is established, the interface starts to propagate in a much slower
time scale ({\it Motion of interface}). The interface obeys a
certain law of motion, which is to be investigated.

A first step to understand this motion is to use (inner and outer)
formal asymptotic expansions of $u^\ep$. This was performed in the
pioneering work of Allen and Cahn \cite{AC} and, slightly later,
in Kawasaki and Ohta \cite{KO}, who revealed that the interface
motion involves curvature effects. Using such arguments, one
discovers that the sharp interface limit of $({\rm P}^\ep)$ obeys
the following law of motion:
\[
 ({\rm P}^0)\quad\begin{cases}
 \, V_{n}=-(N-1)\kappa + c_0 \displaystyle \int _{\am}^{\ap}
 g(x,t,r)\,dr \quad \text { on } \Gamma_t \\
 \, \Gamma_t\big|_{t=0}=\Gamma_0\,,
\end{cases}
\]
where $\Gamma_t$ denotes the limit sharp interface at time $t\geq
0$, $V_n$ is the normal velocity of $\Gamma _t$ in the exterior
direction, $\kappa$ the mean curvature at each point of
$\Gamma_t$, $c_0$ a constant determined straightforwardly from $f$
via
\[
c_0:=\Big [\sqrt 2\int_{\am}^{\ap} (W(s) - W(\am))^{1/2} ds \Big ]
^{-1},
\]
where $W(s):=-\int_a^s f(r) dr$. As long as the solution
$\Gamma_t$ of $({\rm P}^0)$ exists, we denote by $\Omega ^- _t$
the region enclosed by $\Gamma_t$, and by $\Omega ^+ _t$ the
region enclosed between $\partial \Omega$ and $\Gamma _t$. Also we
define a step function $\tilde u(x,t)$ by
\begin{equation}
\tilde u(x,t):=\begin{cases}\label{def:tilde-u}
\, \am &\text{in } \Omega ^-_t\\
\, \ap &\text{in } \Omega ^+_t\,,
\end{cases}
\end{equation}
to which $u^\ep$ is formally supposed to converge as $\ep\to 0$.
As regards $({\rm RD}^\ep)$, the limit problem is found to be
\[
({\rm RD}^0)\quad
\begin{cases}
\,V_{n}=-(N-1)\kappa - c_0 \displaystyle \int_{\am}^{\ap}
f_1(r,v)\,dr \quad &\textrm{on}\ \ \Gamma_t
\vspace{2pt}\\
\,\tilde{v}_t=D \Delta \tilde{v} + h(\tilde{u},\tilde{v})
\quad &\textrm{in}\ \ \Omega\times(0,\infty)\vspace{2pt} \\
\,\Gamma_t\big|_{t=0}=\Gamma_0 \vspace{2pt}\\
\,\displaystyle{\frac{\partial \tilde{v}}{\partial \nu}} = 0
\quad&\textrm{on}\ \ \partial \Omega \times (0,\infty)\vspace{2pt}
\\
\,\tilde v(x,0)=v_0(x) \quad&\textrm{in}\ \ \Omega\,,
\end{cases}
\]
where $\tilde{u}$ is the step function defined in
\eqref{def:tilde-u}. This is a system consisting of an equation of
surface motion and a parabolic partial differential equation.
Since $\tilde{u}$ is determined straightforwardly from $\Gamma_t$,
in what follows, by a solution of $({\rm RD}^0)$ we mean a pair
$(\Gamma,\tilde{v})$.

\begin{remark}[Local smooth solutions for the limit problems]
\label{rm:local-smooth} Under our assumptions, there exists
$T^{max}>0$ such that $({\rm P}^0)$, respectively $({\rm RD}^0)$,
possesses a unique smooth solution $\Gamma=\cup_{0\leq
t<T^{max}}(\Gamma _t\times \{t\})$, resp. $(\Gamma, \tilde
v)=(\cup _{0\leq t<T^{max}}(\Gamma _t\times \{t\}), \tilde v)$.
For more details we refer to \cite{A-Hil-Mat} and the references
therein, in particular \cite{CR}, \cite{Che-xy}, \cite{C2}. In the
sequel we select any $0<T<T^{max}$ and work on $[0,T]$.
\end{remark}

Numerous efforts have been made to rigorously prove the
convergence of $({\rm P}^\ep)$ and $({\rm RD} ^\ep)$ to $({\rm
P}^0)$ and $({\rm RD}^0)$, respectively. Concerning the Allen-Cahn
equation, let us mention the work of de Mottoni and Schatzman
\cite{MS1} (generation of interface via sub- and super-solutions)
and \cite{MS2} (motion of interface via construction of and
linearization around an ansatz) or that of Bronsard and Kohn
\cite{BK} (motion of interface via $\Gamma$-convergence). Chen
\cite{C1,C2} has established an $\mathcal O(\ep |\ln \ep|)$ error
estimate between the location of the actual transition later and
the limit interface, both for scalar equations and systems for
rather general initial data. More recently, in \cite{A-Hil-Mat},
the present authors improved this estimate to $\mathcal O(\ep)$.
More precisely, they show that the solution $u^\ep$ develops a
steep transition layer within the time scale of $\mathcal
O(\ep^2|\ln\ep|)$, and that the layer obeys the law of motion that
coincides with the formal asymptotic limit $({\rm P}^0)$ or $({\rm
RD}^0)$ within an error margin of $\mathcal O(\ep)$. Let us also
mention that there are results of much finer error estimates in
the literature (see for instance \cite{Bel-Pao}), but those
results are concerned with very specific initial data which
already have nice transition layers consistent with formal
asymptotics (hence dependent on $\ep$).

As mentioned before, in most of the aforementioned works,
approximate solutions or sub- and super-solutions are constructed
by roughly following the formal expansions. Our goal is to
investigate the actual validity of such expansions for solutions
$u^\ep$ with rather general initial data.

\begin{remark}[Viscosity framework]
Since the limit problem may develop singularities in finite time,
the {\it classical framework} does not always allow to study the
singular limit procedure for all $t\geq 0$. Nevertheless --- as
far as the Allen-Cahn equation is concerned --- following
\cite{ES1}, \cite{CGG} one can define a limit problem for all
$t\geq 0$ that generalizes $({\rm P}^0)$ in the framework of {\it
viscosity solutions}. In this setting we refer to \cite{ESS}
(convergence of Allen-Cahn equation with prepared initial data to
generalized motion by mean curvature), \cite{BBS}, \cite{BSS}
(generalizations), \cite{S2, S3}, \cite{BaS}, \cite{BD} (not
well-prepared initial data), \cite{A-Dro-Mat} (fine convergence
rate).
\end{remark}

%%%%%%%%%%%%%%%%%%%%%%%%%%%%%
%%%%%%%%%%%%%%%%%%%%%%%%%%%%%
\section{Main results}\label{s:results}

We start by giving an outline of the formal asymptotic expansions
mentioned before. See \cite[Section 2]{A-Hil-Mat} for more
details. Let $u^\ep$ be the solution of $({\rm P}^\ep)$, and
$\Gamma=\cup _{0\leq t \leq T} (\Gamma _t \times \{t\})$ be the
solution of the limit geometric motion problem $({\rm P}^0)$. We
define the {\it signed distance function} to $\Gamma$ by
\begin{equation}\label{eq:dist}
 d (x,t):=
 \begin{cases}
 -&\hspace{-10pt} \mbox{dist}(x,\Gamma _t) \quad \text{for }
 x\in\Omega _t^- \\
 &\hspace{-10pt} \mbox{dist}(x,\Gamma _t)\quad\text{for }x\in\Omega
 _t^+\,.
 \end{cases}
\end{equation}
Then, near $\Gamma$, we make a formal inner expansion of the form
\begin{equation}\label{expansion}
 u^\ep(x,t)=U_0\left(x,t,\frac{ d(x,t)}{\ep}\right)+\ep
 U_1\left(x,t,\frac{ d(x,t)}{\ep}\right)+\cdots\,.
\end{equation}
Some normalization conditions and matching conditions (with the
outer expansion) are also imposed. By plugging the expansion
\eqref{expansion} into $({\rm P}^\ep)$, we discover that
$U_0(x,t,z)=U_0(z)$, where $U_0(z)$ is the unique solution (whose
existence is guaranteed by the integral condition \eqref{int-f})
of the stationary problem
\begin{equation}\label{eq-phi}
 \left\{\begin{array}{ll}
 {U_0} '' +f(U_0)=0 \vspace{2pt}\\
 U_0(-\infty)= \am\,,\quad U_0(0)=a\,,\quad U_0(\infty)= \ap\,.
 \end{array} \right.
\end{equation}
This solution represents the first approximation of the profile of
a transition layer around the interface observed in the stretched
coordinates. Next the solvability condition for the equation
involving $U_1$ provides the law of motion $({\rm P}^0)$ for the
limit interface $\Gamma$, which, in turn, determines the term
$U_1$.

It is then natural to wonder if the ansatz
\[
u^\ep (x,t)=U_0\left(\frac{d(x,t)}{\ep}\right)+\ep
U_1\left(x,t,\frac{d(x,t)}{\ep}\right)+\cdots
\]
is really a good approximation of the profile of the solution
$u^\ep$. Note that the convergence results mentioned in subsection
\ref{ss:singular} do not answer this question; indeed those
results simply show that the level surface of the solution $u^\ep$
\begin{equation}\label{level-sets}
\Gamma _t ^\ep:=\{x\in \Omega:\,u^\ep(x,t)=a\}
\end{equation}
converges to the sharp interface $(\Gamma _t)_{0\leq t \leq T}$,
which is a solution of $({\rm P}^0)$, and that
\begin{equation}\label{ue-to-u}
 \lim_{\ep\to 0}u^\ep(x,t)=
 \begin{cases}
 \alpha^- \quad \text{for } x\in\Omega_t^-\\
 \alpha^+ \quad\text{for } x\in\Omega_t^+\,,
 \end{cases}
\end{equation}
without clarifying the validity of \eqref{expansion}. Our main
result Theorem \ref{th:validity-allen} below provides a first
answer in this direction.
%%It asserts that, as $\ep \to 0$,
%%$(i)$ the level set $\Gamma _t ^\ep$ turns out to be a smooth
%%hypersurface, so that we can define $d^\ep$ the associated signed
%%distance function to $\Gamma ^\ep$; $(ii)$---$(iii)$ the term
%%$U_0(\frac{d^\ep(x,t)}\ep)=U_0(\frac{d(x,t)+\mathcal O(\ep)}\ep)$
%%is actually a good approximation of $u^\ep(x,t)$.
In the sequel we define
\begin{equation}\label{time}
t^\ep:=f'(a)^{-1}\ep ^2|\ln \ep |\,,
\end{equation}
which is the time needed for the transition layer of $u^\ep$ to
become fully well-developed (see Lemma \ref{lem:AHM}). We define
the {\it signed distance function associated with $\Gamma^\ep$} by
\begin{equation}\label{eq:dist-e}
 d^\ep (x,t):=
 \begin{cases}
 -&\hspace{-10pt} \mbox{dist}(x,\Gamma^\ep _t) \quad \text{if }
 u^\ep(x,t)<a \\
 &\hspace{-10pt} \mbox{dist}(x,\Gamma^\ep _t)\quad\text{if }
 u^\ep(x,t)>a\,.
 \end{cases}
\end{equation}
Note that this definition of $d^\ep$ is consistent with that of
$d$ in \eqref{eq:dist} in view of \eqref{ue-to-u}.

\begin{theorem}[Validity for Allen-Cahn]\label{th:validity-allen}
Let the assumptions of subsection \ref{ss:assumptions} hold (in
particular the initial condition $u_0$ is rather generic). Let
$u^\ep$ be the smooth solution of Allen-Cahn equation $({\rm
P}^\ep)$. Fix $\mu >1$. Then the following hold.
\begin{enumerate}
\item [{\rm (i)}] If $\ep >0$ is small enough then, for any
$t\in[\mu t^\ep,T]$, the level set $\Gamma _t ^\ep$ is a smooth
hypersurface and can be expressed as a graph over $\Gamma _t$.
\item[{\rm (ii)}]
\begin{equation}\label{lim-d-epsilon}
\lim_{\ep \to 0}\; \sup_{\mu t^\ep \leq t\leq T,\, x\in \bar
\Omega} \left |u^\ep(x,t)-U_0\left(\frac
{d^\ep(x,t)}{\ep}\right)\right|=0\,,
\end{equation}
where $d^\ep$ denotes the signed distance function associated with
$\Gamma^\ep$. \item [{\rm (iii)}] There exists a family of
functions
\begin{equation}\label{theta}
  \theta ^\ep:\cup _{0\leq t \leq T} (\Gamma _t \times \{t\})\to
  \R \quad (0<\ep<<1)
\end{equation}
whose $L^\infty$-norms remain bounded as $\ep \to 0$, such that
\begin{equation}\label{lim}
\lim_{\ep \to 0}\; \sup_{\mu t^\ep \leq t\leq T,\, x\in \bar
\Omega} \left |u^\ep(x,t)-U_0\left(\frac {d(x,t)-\ep \theta
^\ep(p(x,t),t)}{\ep}\right)\right|=0\,,
\end{equation}
where $d$ denotes the signed distance function associated with
$\Gamma$ and $p(x,t)$ denotes a point on $\Gamma _t$ such that
$\mbox{dist} (x,\Gamma _t)=\Vert x-p(x,t)\Vert$.
\end{enumerate}
\end{theorem}

Note that $p(x,t)$ is an orthogonal projection of the point $x$
onto $\Gamma _t$, which is uniquely defined in a small tubular
neighborhood of $\Gamma _t$ since $\Gamma_t$ is a smooth solution
of $({\rm P}^0)$. Note also that the presence of the perturbations
$-\ep \theta ^\ep(p(x,t),t)$ cannot be avoided since it reflects
the small difference between $d(x,t)$ and $d^\ep (x,t)$.

Let us mention that the validity of higher order terms of the
formal expansions (for generic solutions) is still unknown.

Our next theorem provides similar estimates for the
reaction-diffusion systems.

\begin{theorem}[Validity for the reaction-diffusion system]
\label{th:validity-nagumo} Let the assumptions of subsection
\ref{ss:assumptions} hold, and let $(u^\ep,v^\ep)$ be the smooth
solution of the reaction-diffusion system $({\rm RD}^\ep)$. Fix
$\mu >1$. Then the same conclusions as in Theorem
\ref{th:validity-allen} hold, with $d$ being the signed distance
function associated with $(\Gamma,\tilde v)$, which is the smooth
solution of $({\rm RD}^0)$ on $[0,T]$.
\end{theorem}

%%%%%%%%%%%%%%%%%%%%%%%%%%%
%%%%%%%%%%%%%%%%%%%%%%%%%%%
\section{Proof of the main results}\label{s:proof}

The proof of Theorem \ref{th:validity-allen} relies on the
following two results:
\begin{quote}
\begin{itemize}
\item [(a)] {\it the level set $\Gamma ^\ep _t$ is approximated by
the interface $\Gamma _t$ by order $\mathcal O(\ep)$
(\cite{A-Hil-Mat}, see subsection \ref{ss:thickness} of the
present paper),} \item[(b)] {\it any eternal solution that lies
between two planar waves is actually a planar wave
(\cite{Ber-Ham}, see subsection \ref{ss:eternal} of the present
paper),}
\end{itemize}
\end{quote}
combined with a rescaling argument.
%%Note that a related argument can be found in \cite{Mat-Nar}
%%where perturbations of Allen-Cahn planar fronts are studied.

%%%%%%%%%%%%%%%%%%%%%%%%%%%%%
\subsection{Thickness of the layers: the refined $\mathcal O(\ep)$
estimate}\label{ss:thickness}

We  quote a result which is valid for both $({\rm P}^\ep)$ and
$({\rm RD}^\ep)$.

\begin{lemma}[{\cite[Theorem 1.3 and Theorem 1.11]{A-Hil-Mat}}]
\label{lem:AHM} Let $\eta$ be an arbitrary constant satisfying $0<
\eta <\min(a-\am,\ap -a)$. Then there exist positive constants
$\ep _0 $ and $C_0$ such that, for all $\,\ep \in (0,\ep _0)$ and
for all $\,f'(a) ^{-1} \ep ^2 |\ln \ep|=t^\ep \leq t \leq T$, we
have
\begin{equation}\label{thickness}
|u^\ep(x,t)-\apm|\leq \eta\quad \text{ if }\, x\in \Omega _t ^\pm
\setminus \mathcal N_{C_0 \ep}(\Gamma _t)\,,
\end{equation}
where $\mathcal N _r(\Gamma _t):\, =\{x\in \Omega:\, dist(x,\Gamma
_t)<r\}$ denotes the $r$-neighborhood of $\Gamma _t$. This implies
in particular that $\Gamma _t ^\ep \subset \mathcal N_{C_0
\ep}(\Gamma _t)$ for all $t^\ep \leq t\leq T$, hence
\begin{equation}\label{d-de}
|d^\ep(x,t)-d(x,t)|\leq C_0\ep \quad \hbox{for all}\
(x,t)\in\overline\Omega \times[t^\ep,T]\,,\ 0<\ep<<1\,.
\end{equation}
\end{lemma}

%%%%%%%%%%%%%%%%%%%%%%%%%%%
\subsection{Eternal solutions and planar waves}\label{ss:eternal}

We recall that a solution of an evolution equation is called {\it
eternal} (or an {\it entire} solution) if it is defined for all
positive and negative time. We follow this terminology to refer to
a solution $w(z,\tau)$ of
\begin{equation}\label{eq-eternal}
w_{\tau} = \Delta _z w + f(w)\,, \quad z\in\R^N,\,\tau \in\R\,.
\end{equation}

Stationary solutions and travelling waves are examples of eternal
solutions. Crucial to our analysis is a recent result of
Berestycki and Hamel \cite{Ber-Ham} asserting that \lq\lq any
planar-like eternal solution is actually a planar wave". More
precisely, the following holds (for $z\in \R ^N$ we write
$z=(z^{(1)},\cdots,z^{(N)})$).

\begin{lemma}[{\cite[Theorem 3.1]{Ber-Ham}}]\label{lem:Ber-Ham}
Let $w(z,\tau)$ be an eternal solution of \eqref{eq-eternal}
satisfying
\begin{equation}\label{conditions}
\liminf_{z^{(N)}\to \infty}\; \inf_{z'\in\R^{N-1}} w(z,\tau) >
a\,, \quad \limsup_{z^{(N)}\to-\infty}\; \sup_{z'\in\R^{N-1}}
w(z,\tau) < a\,,
\end{equation}
where $z':=(z^{(1)},\cdots,z^{(N-1)})$. Then there exists  a
constant $z^*\in\R$ such that
\[
w(z,\tau) = U_0(z^{(N)}-z^*)\,,  \quad z\in\R^N\,,\,\tau\in\R\,.
\]
\end{lemma}

%%%%%%%%%%%%%%%%%%%%%%%%%
\subsection{Proof of (ii) in Theorem \ref{th:validity-allen}}
\label{SS:proof1}

In what follows we fix $\mu>1$ and an arbitrary constant $T_1$
with $T<T_1<T^{max}$ (see Remark \ref{rm:local-smooth}). Obviously
the conclusion of Lemma \ref{lem:AHM} remains valid if $T$ is
replaced by $T_1$. Assume by contradiction that
\eqref{lim-d-epsilon} does not hold. Then there is $\eta >0$ and
sequences $\ep _k \downarrow 0$, $t_k\in[\mu t^{\ep _k},T]$, $x_k
\in \bar \Omega$ ($k=1,2,...$) such that
\begin{equation}\label{eta}
  \left |u^{\ep _k}(x_k,t_k)-U_0\left (\frac{d^{\ep_k}(x_k,t_k)}
  {\ep_k}\right)\right|\geq 2\eta\,.
\end{equation}
In view of \eqref{thickness}--\eqref{d-de} and $U_0(\pm
\infty)=\apm$, for \eqref{eta} to hold it is necessary to have
\begin{equation}\label{near}
  d(x_k,t_k)=\mathcal O (\ep _k)\,,\quad \text{ as }\, k\to
  \infty\,.
\end{equation}

If $u^{\ep _k}(x_k,t_k)=a$, then this would mean that $x_k\in
\Gamma^{\ep_k}_{t_k}$, in which case the left-hand side of
\eqref{eta} would be $0$ (since $U_0(0)=a$), which is impossible.
Hence $u^{\ep _k}(x_k,t_k)\ne a$. By extracting a subsequence if
necessary, we may assume without loss of generality that $u^{\ep
_k}(x_k,t_k)- a$ has a constant sign for $k=0,1,2,\ldots$. Since
the sign of this quantity is irrelevant in the later argument, in
what follows we assume that
\begin{equation}\label{u<0}
  u^{\ep _k}(x_k,t_k)> a \quad (k=0,1,2,\ldots),
\end{equation}
which then implies that
\[
  d^{\ep_k}(x_k,t_k)> 0 \quad (k=0,1,2,\ldots).
\]
Since the sequence $(x_k)$ remains close to $\Gamma_{t_k}$ by
\eqref{near}, and since $\Gamma_{t_k}$ is uniformly smooth for
$k=0,1,2,\ldots$, each $x_k$ has a unique  orthogonal projection
$p(x_k,t_k) \in \Gamma _{t_k}$.  Let $y_k$ be a point on
$\Gamma^{\ep_k}_{t_k}$ that has the smallest distance from $x_k$.
If such a point is not unique, we choose one such point
arbitrarily. Then we have
\begin{equation}\label{u=0}
  u^{\ep_k}(y_k,t_k)=a \quad (k=0,1,2,\ldots),
\end{equation}
\begin{equation}\label{d=x-y}
  d^{\ep_k}(x_k,t_k)=\Vert x_k - y_k\Vert\,,
\end{equation}
\begin{equation}\label{circle}
  u^{\ep_k}(x,t_k)>a \quad\  \hbox{if}\ \ \Vert x-x_k\Vert
  < \Vert y_k - x_k\Vert\,,
\end{equation}
\[
x_k-p_k\perp \Gamma_{t_k} \quad\hbox{at}\ \ p_k\in\Gamma_{t_k}\,,
\]
where $p_k:=p(x_k,t_k)$. Furthermore, \eqref{near} and
\eqref{d-de} imply
\begin{equation}\label{near2}
  \Vert x_k-p_k\Vert = {\mathcal O}(\ep_k)\,, \ \
  \Vert y_k-p_k\Vert = {\mathcal O}(\ep_k) \quad
  (k=0,1,2,\ldots).
\end{equation}

We now rescale the solution $u^\ep$ around $(p_k,t_k)$ and define
\begin{equation}\label{rescaling}
  w^k(z,\tau):=u^{\ep _k}(p_k+\ep _k \mathcal R _k z,t_k+\ep _k ^2
  \tau)\,,
\end{equation}
where $\mathcal R _k$ is a matrix in $SO(N,\R)$ that rotates the
$z^{(N)}$ axis onto the normal at $p_k\in \Gamma_{t_k}$, that is
\[
\mathcal  R_k:(0,\dots,0,1)^T\mapsto n(p_k,t_k)\,,
\]
where $(\ )^T$ denotes a transposed vector and $n(p,t)$ the
outward normal unit vector at $p\in \Gamma_{t}$. Since $\Gamma_t$
(hence the points $p_k$) is uniformly separated from
$\partial\Omega$ by some positive distance, there exists $c>0$
such that $w^k$ is defined (at least) on the box
\[
 B^k:=\left \{(z,\tau)\in \R^N\times \R \,:\, \Vert z\Vert \leq
 \frac c{\ep _k}, \ \  -(\mu -1)f'(a)^{-1}|\ln \ep
_k|\leq \tau \leq \frac{T_1-T}{\ep _k ^2}\right\}\,.
\]
Since $u^\ep$ satisfies $({\rm P}^\ep)$, we see that $w^k$
satisfies
\begin{equation}\label{edp-wk}
  w^k _\tau=\Delta _z w^k+f(w^k)-\ep_k g^{\ep _k}(p_k+\ep _k
  \mathcal R _k z, t_k+\ep _k ^2 \tau,w^k)\quad \text{ in } B^k\,.
\end{equation}
Moreover, if $(z,\tau)\in B^k$ then $t^{\ep_k}\leq t_k+\ep _k ^2
\tau \leq T_1$. Therefore \eqref{thickness} implies
%%(recall that $x_k-y_k$ points towards$\Omega _{t_k}^+$)
\begin{equation}\label{front-like0}
  \left\{ \begin{array}{lll}
  d(p_k+\ep _k \mathcal R _k z,t_k+\ep _k^2\tau) \leq -C_0\ep _k \ \ & \Rightarrow
  \ & w^k(z,\tau) \leq \am+\eta\,,\vspace{4pt}\\
  d(p_k+\ep _k \mathcal R _k z,t_k+\ep _k^2\tau) \geq C_0\ep _k \ \ & \Rightarrow
  \ & w^k(z,\tau) \geq \ap-\eta\,,
  \end{array}\right.
\end{equation}
so long as $(z,\tau) \in B^k$. Since the rotation by $\mathcal R
_k$ of the $z^{(N)}$ axis is normal to $\Gamma_{t_k}$ at
$p_k:=p(x_k,t_k)$, and since the curvature of $\Gamma_t$ is
uniformly bounded for $0\leq t\leq T$, we see from
\eqref{front-like0} that there exists a constant $C>0$ such that
\begin{equation}\label{front-like}
  z^{(N)}\leq -C \Rightarrow w^k(z,\tau)\leq \am+\eta\,,\quad
   z^{(N)}\geq C \Rightarrow w^k(z,\tau)\geq \ap-\eta\,,
\end{equation}
so long as $(z,\tau) \in B^k$ and $\Vert z\Vert\leq
\sqrt{\ep_k^{-1}}$.

Now, since $w^k$ solves \eqref{edp-wk}, the uniform (w.r.t. $k\geq
0$) boundedness of $w^k$ and standard parabolic estimates, along
with the derivative bounds on $g^\ep$, imply that $w^k$ is
uniformly bounded in $C_{loc}^{2+\gamma,1+\frac \gamma 2}(B^1)$.
We can therefore extract from $(w^k)$ a subsequence that converges
to some $w$ in $C_{loc}^{2,1}(B^1)$. By repeating this on all
$B^k$, we can find a subsequence of $(w^k)$ that converges to some
$w$ in $C^{2,1}_{loc}(\R^N\times \R)$ (note that $\cup_{k\geq 0}
B^k=\R^N \times \R$). Passing to the limit in \eqref{edp-wk}
yields
\[
 w_\tau=\Delta _zw +f(w)\quad \text{ on } \R^N\times \R\,.
\]

Hence we have constructed an eternal solution $w(z,\tau)$ which
--- in view of \eqref{front-like}--- satisfies
\eqref{conditions}. Lemma \ref{lem:Ber-Ham} then implies that
\begin{equation}\label{w=U0}
w(z,\tau)=U_0(z^{(N)}-z^*)
\end{equation}
for some $z^*\in\R$.

Now we define sequences of points $(z_k),\,(\tilde z_k)$ by
\[
 z_k:=\frac 1 {\ep _k} \mathcal R _k ^{-1}(x_k-p_k)\,,\quad
 \tilde z_k:=\frac 1 {\ep _k} \mathcal R _k ^{-1}(y_k-p_k)\,.
\]
By \eqref{near2}, these sequences are bounded, so we may assume
without loss of generality that they converge:
\[
  z_k\to z_\infty\,, \quad \tilde z_k\to \tilde z_\infty\,, \quad
  \hbox{as}\ \ k\to\infty\,.
\]
By the definition of the $z$ coordinates, $z_\infty$ must lie on
the $z^{(N)}$ axis, that is,
\[
  z_\infty=(0,\dots,0,z^{(N)}_\infty)^T\,.
\]
It follows from \eqref{u=0} and \eqref{circle} that
\begin{equation}\label{circle-w}
  w(\tilde z_\infty,0)=a\,,\quad\ \  w(z,0)\geq a \ \ \hbox{if}\ \
  \Vert z-z_\infty\Vert \leq     \Vert\tilde z_\infty-z_\infty\Vert\,.
\end{equation}
Note that by \eqref{w=U0}, the level set $w(z,0)=a$ coincides with
the hyperplane $z^{(N)}=z^*$, and recall that ${U_0}'>0$.
Therefore, in view of \eqref{w=U0} and \eqref{circle-w}, we have
either $\tilde z_\infty=z_\infty$, or that the ball of radius
$\Vert\tilde z_\infty-z_\infty\Vert$ centered at $z_\infty$ is
tangential to the hyperplane $z^{(N)}=z^*$ at $\tilde z_\infty$.
This implies that $\tilde z_\infty$, as well as $z_\infty$, must
also lie on the $z^{(N)}$ axis. Therefore
\[
\tilde z_\infty=(0,\dots,0,z^*)^T\,,
\]
and the inequality $w(z_\infty,0)\geq a$ implies that
$z_\infty^{(N)}\geq z^*$. On the other hand \eqref{d=x-y} implies
$d^{\ep_k}(x_k,t_k)/\ep_k=\Vert x_k-y_k\Vert/\ep_k=\Vert
z_k-\tilde z_k\Vert\to \Vert z_\infty-\tilde z_\infty \Vert
=z^{(N)}_\infty-z^*$. The assumption \eqref{eta} then yields
\begin{eqnarray*}
0&=&\Vert w(z_\infty,0)-U_0(z_\infty^{(N)}-z^*)\Vert\\
&=&\Vert \lim_{k\to\infty}
u^{\ep_k}(x_k,t_k)-U_0\Big(\lim_{k\to\infty}\frac{d^{\ep_k}(x_k,t_k)}
{\ep_k}\Big)\Vert\\
& \geq& 2\eta\,.
\end{eqnarray*}
This contradiction proves statement (ii) of Theorem \ref{lim}.

%%%%%%%%%%%%%%%%%%%%%
\subsection{Proof of (i) and (iii) in Theorem \ref{th:validity-allen}}

The proof of (i) below uses an argument similar to the proof of
Corollary 4.8 in \cite{Mat-Nar}. Fix $\mu >1$. For a given $\eta
\in (0,\min(a-\am,\ap -a))$ define $\ep _0 >0$ and $C_0 >0$ as in
Lemma \ref{lem:AHM}. Then we claim that
\begin{equation}\label{claim}
  \liminf _{\ep \to 0} \inf _{x\in \mathcal N_{C_0 \ep}(\Gamma
  _t),\, \mu t^\ep \leq t \leq T} \nabla u^\ep (x,t) \cdot
  n(p(x,t),t)>0\,,
\end{equation}
where $n(p,t)$ denotes the outward unit normal vector at $p\in
\Gamma _t$. Indeed, assume by contradiction that there exist
sequences $\ep _k \downarrow 0$, $t_k\in[\mu t^{\ep _k},T]$, $x_k
\in \mathcal N_{C_0 \ep _k}(\Gamma _{t_k})$ ($k=1,2,...$) such
that
\[
  \nabla u^{\ep _k} (x_k,t_k) \cdot n(p_k,t_k) \leq 0\,,
\]
where $p_k=p(x_k,t_k)$. By rescaling around $(p_k,t_k)$ and using
arguments similar to those in the proof of (ii), one can find a
point $z_\infty$ with $|z_\infty ^{(N)} |\leq C_0$ such that
\[
  {U_0}'(z_\infty ^{(N)})\leq 0\,,
\]
which contradicts the fact that ${U_0}'>0$ and establishes
\eqref{claim}. Since, in view of Lemma \ref{lem:AHM}, $\Gamma _t
^\ep \subset \mathcal N _{C_0 \ep}(\Gamma _t)$, the estimate
\eqref{claim} implies that $\nabla u^\ep (x,t) \neq 0$ for all
$x\in \Gamma _t ^\ep$; hence by the implicit function theorem,
$\Gamma _t ^\ep$ is a smooth hypersurface in a neighborhood of any
point on it. The fact that $\Gamma _t ^\ep$ can be expressed as a
graph over $\Gamma _t$ also follows from \eqref{claim}. This
proves statement (i) of Theorem \ref{th:validity-allen}.

Finally, statement (iii) follows immediately from statements (i),
(ii) and \eqref{d-de}. This completes the proof of Theorem
\ref{th:validity-allen}.

%%%%%%%%%%%%%%%%%%%%%%%%%
\subsection{Proof of Theorem \ref{th:validity-nagumo}}

As shown in \cite[Section 7]{A-Hil-Mat}, the behavior of $u^\ep$
in the system $({\rm RD}^\ep)$ can be treated as a special case of
$({\rm P}^\ep)$, by regarding $v^\ep$ as a given function and
using a contraction mapping theorem. Thus the conclusion of
Theorem \ref{th:validity-nagumo} follows directly from Theorem
\ref{th:validity-allen}.

%%For acknowledgements section, please don't number the section,
%%please begin it with \section*{Acknowledgements}
%%\section*{Acknowledgments} The authors would like to thank the anonymous refer%%ee for his valuable comments.

%%%%%%%%%%%%%%%%%%%%%%%%%%%%%%
%%%%%%%%%%%%%%%%%%%%%%%%%%%%%%

\end{document}